\documentclass[12pt,intlim,righttag]{article}
\usepackage{amsmath,amsthm}
\usepackage{amssymb}
\usepackage{amsfonts}

\newcommand{\la}{\langle}
\newcommand{\ra}{\rangle}

\newcommand{\vp}{\varphi}

\newcommand{\td}{\tilde}

\newcommand{\be}{\begin}
\newcommand{\ee}{\end}

\newcommand\beq{\begin{equation}}
\newcommand\eeq{\end{equation}}

\newcommand{\bea}{\begin{eqnarray}}
\newcommand{\eea}{\end{eqnarray}}
\newcommand{\beaa}{\begin{eqnarray*}}
\newcommand{\eeaa}{\end{eqnarray*}}

\theoremstyle{Theorem}

\theoremstyle{corollary}

\theoremstyle{remark}

\theoremstyle{definition}

\def\a{\alpha}
\def\b{\beta}

\def\d{\delta}

\begin{document}
\title{On Martingale Transformations of Multidimensional Brownian Motion
 }

\author{M. Mania$^{1)}$ and R. Tevzadze$^{2)}$}

\date{~}
\maketitle

\begin{center}
$^{1)}$ A. Razmadze Mathematical Institute of Tbilisi State University  and
Georgian-American University, Tbilisi, Georgia,
\newline(e-mail: misha.mania@gmail.com)
\\
$^{2)}$ Georgian-American University and
Institute of Cybernetics of Georgian Technical Univercity,  Tbilisi, Georgia
\newline(e-mail: rtevzadze@gmail.com)
\end{center}

\begin{abstract}
{\bf Abstract.} We describe the class of functions $f: R^n\to R^m$ which transform
a vector Brownian Motion into a martingale and use this description to give martingale
characterization of the general measurable solution of the multidimensional Cauchy functional equation.

\end{abstract}



\section{Introduction}

It is well known (see, e.g., \cite{W}, \cite{CI}, \cite{Ch}, \cite{MT}) that if $f=(f(x), x\in R)$ is a function of one variable and $W$ is linear Brownian Motion  then
the transformed process $f(W_t)$ is
a continuous (or right-continuous) martingale if and only if $f$ is an affine function. In multidimensional case this result is no longer true.
A simple countre-example gives the non-linear function $f(x,y)=x^2-y^2$  for which the transformed process $f(W^1_t, W_t^2)$ of two independent
Brownian Motions $W^1$ and $W^2$ is a continuous martingale.
Our goal is  to give  sufficient (and necessary) conditions in multidimensional case, when martingale  function is affine or almost affine,
i.e., it coincides with an affine function almost everywhere with respect to the Lebesgue measure.

Let $W=(W_t, t\ge 0)$ be a  $n$-dimensional standard Brownian Motion defined on a  probability space  $(\Omega, \cal F, P)$ with filtration
  $F=({{\cal F}}_t,t\ge0)$  satisfying  the usual conditions of right-continuity and completeness.

 In \cite{BCD} the functions $f:R^n\to R^m$ was studied such that $f(W)$ is  Brownian path preserving, i.e., is a standard Brownian Motion up to a random time change.
 It follows from their results that, if $m=1$ and $f$ is continuous, then  the process $f(W)$ is Brownian path preserving if and only if $f$ is harmonic. It was shown in \cite{T} that
 if $f(W)$ is again a Brownian motion (without allowing time change) with respect to the same filtration,  then $f$ is an affine function.  In Theorem 1  we consider the case
 when the transformed  process $f(W)$ is a general martingale, without assuming the continuity of paths, but impose an additional condition, which is satisfied in the case of Brownian motion and which
 guaranties the almost linearity of the function $f$ and the linearity if we additionally  assume the continuity of $f$.

 In section 3 we apply these results to give an equivalent martingale characterization of the general measurable solution of multidimensional Cauchy's functional equation (see, e. g., \cite {AD}, \cite{SK} and  \cite{J}, \cite{DB}
 for almost additive version). In Theorem 3 we show
 that if $f=(f(x), x\in R^n)$ is a measurable function  satisfying the  Cauchy functional equation
$$
f(x+y)=f(x)+f(y),
$$
for almost all $(x,y)$ in the sense of the Lebesgue measure on $R^{2n}$, then the transformed process $f(W_t)$ is
a  martingale which satisfies condition C) of Theorem 1 and, hence  is almost affine function.

\section{Martingale functions of Brownian Motion}

Let $M$ be a martingale with respect to the filtration $F^W=({\cal F}^W_t, t\ge0)$  generated
by the Brownian Motion $W$ and denote by $\tilde M$ the continuous modification of $M$. Since
almost all paths of $\tilde M$  are continuous, $\tilde M$ is locally square integrable martingale and the square characteristic $\la {\tilde M}\ra$ of
$\tilde M$ exists.  We shall call $\la{\tilde M}\ra$ the square characteristic  of $M$ also, i.e.,   $\la M\ra\equiv\la {\tilde M}\ra$.

{\bf Theorem 1.} Let $f(x)=(f_1(x),...,f_m(x)), x\in R^n$ be a measurable function, such that $f(W_t)=(f_1(W_t),...,f_m(W_t))$ is a martingale
satisfying condition:
$$
C)\;\text{the process}\;\; \la f_j(W)\ra_t-C_j t\;\;\text{is a non-decreasing for some}\;\;C_j,\;j=1,...,n.
$$

Then $f(x)$  must be an affine function for almost all $x$ in the sense of the Lebesgue measure on $R^{n}$.
\begin{proof}
Let
\beaa
g\left(t,x\right)=E\left(f\left(W_{T}\right)|W_{t}=x\right).
\eeaa
It is evident that $g\left(t,x\right)$ satisfies the generalized heat equation
\begin{equation}\label{gen}
\int_0^T\int_{R^n} g(s,y)\big(\frac{\partial\varphi}{\partial t}(s,y)-\frac{1}{2}\Delta\varphi(s,y)\big)dyds=0,
\end{equation}
for every infinitely differentiable finite (on $[0,T]\times {R^n}$) function $\varphi$.

By the Markov property of the Brownian motion
\beaa
g\left(t,W_{t}\right)=E\left(f\left(W_{T}\right)|\mathcal{F}_{t}\right)\;\;\;\;\text{a.s.}
\eeaa
and from the martingale property of $f\left(W_{t}\right)$ we have that for all $t\le T$

\beaa
g\left(t,W_{t}\right)=f\left(W_{t}\right) \quad  a.s.
\eeaa

Therefore, for all $t\le T$
\beaa
\int_{R^n}{|g\left(t,x\right)-f\left(x\right)|\frac{1}{\sqrt{2\pi t}}e^{-\frac{|x|^{2}}{2t}}dx}=0
\eeaa
which implies that for all $t\le T$

\beaa
g\left(t,x\right)=f\left(x\right) \quad a.e
\eeaa
with respect to the Lebesgue measure.

Therefore, it follows from (\ref{gen}) that
\begin{equation}\label{gen2}
\int_0^T\int_{R^n} f(y)\big(\frac{\partial\varphi}{\partial t}(s,y)-\frac{1}{2}\Delta\varphi(s,y)\big)dyds=0,
\end{equation}
for every infinitely differentiable finite (on $[0,T]\times R^n$) function $\varphi$.

Taking $\varphi(t,y)=\vp_0(t)\vp_1(y)$ with $\int_0^T\vp_0(t)dt=1$ we obtain from (\ref{gen2}) that
\begin{equation}\label{gen3}
\int_{R^n} f(y)\Delta\varphi_1(y)dy=0,
\end{equation}
for every infinitely differentiable finite on $R^n$ function $\varphi_1$.
From Theorem of \cite{vl} follows that there exists infinite differentiable (even analytic) harmonic function $\td f$ such that $f(x)=\td f(x)$-a.e. .
It is clear that $\td f(W_t)$ is continuous modification of $f(W_t)$. By the Ito formula we get
$$\td f_j(W_t)=\td f(0)+\int_0^t\nabla \td f_j(W_s)dW_s,$$
which means $\la \td f_j(W)\ra_t=\int_0^t\td f_j(W_s)^2ds$.
The condition of Theorem gives $|\nabla \td f_j(W_s)|^2\le C_j^2\;a.e.$ which is the same as $|\nabla \td f_j(x)|\le C_j$. By the mean value theorem we obtain
$$
|\td f_j(x)|\le|\td f_j(0)|+C_j|x|.
$$
Therefore, the Liouville Theorem (\cite{vl}, 290p.) implies that among  harmonic functions in the whole space only affine functions satisfy this condition.
\end{proof}

{\bf Corollary 1.} Let conditions of Theorem 1 are satisfied. If in addition  $f$ is continuous, then it coincides with an affine function.

Sometimes another kind of condition will  be  useful.
Let us introduce the notation
$$W_t^{ij}=\left(W_1(t),...,W_{i-1}(t), \frac{W_{i}(t)+W_{j}(t)}{\sqrt2}, W_{i+1}(t),...,W_n(t)\right).$$

{\bf Theorem 2.} Let $f(x)=(f_1(x),...,f_m(x)), x\in R^n$ be a measurable function, such that  the processes $(f(W_t^{ij}), t\ge0)$ are martingales for each $i\le j$.
Then $f(x)$  is an affine function for almost all $x$ in the sense of the Lebesgue measure on $R^n$.

{\it Proof}. As in Theorem 1 we can prove that $f$ a.e. coincides with an analytic function $\tilde f$ satisfying
$\Delta{\tilde f}(x)+\sqrt2^{1-\d_{ij}}\frac{\partial^2{\tilde f}}{\partial {x_i}\partial{x_j}}(x)=0$
\footnote{$\d_{ij}$  denotes Kroneker's delta} for each
$i,j$.  If we take $i=j$ and subtract equations for different $i$ and $j$ we get $\frac{\partial^2{\tilde f}}{\partial {x_i}^2}(x)=\frac{\partial^2{\tilde f}}{\partial {x_j}^2}(x)$.
Hence
$$\frac{\partial^2{\tilde f}}{\partial {x_i}^2}(x)=\frac1{n+1}(\Delta{\tilde f}(x)+\frac{\partial^2{\tilde f}}{\partial {x_i}^2}(x))=0$$
  and
$$\frac{\partial^2{\tilde f}}{\partial {x_i}\partial{x_j}}(x)=\frac1{\sqrt2}\Delta{\tilde f}(x)+\frac{\partial^2{\tilde f}}{\partial {x_i}\partial{x_j}}(x)=0\;for\;i\neq j.$$
 It follows from these equalities  that  $\frac{\partial^2{\tilde f}}{\partial {x_i}\partial{x_j}}(x)=0$ for any $i,j$, which implies that $\tilde f$ is  an affine function. Hence, $f$ coincides with an affine
 function  almost everywhere with respect to the Lebesgue measure on $R^n$

{\bf Corollary 2.} Let $f=(f(x),x\in R)$ be a function of one variable and $W$ is a linear Brownian Motion. Then

a) If  the process $f(W_t)$ is a martingale, then  the function $f$ coincides with an  affine function almost everywhere  with respect to the Lebesgue measure on $R$.

b) If $f(W_t)$ is a continuous (or right-continuous) martingale then $f$ is  an affine function.

{\it Proof.} a) follows from Theorem 2.
Assertion b) follows from Theorem 2 and the fact that if the process $f(W_t)$ is right-continuous, then the function
$f=(f(x),x\in R)$ will be continuous (see Lemma  A1 from the Appendix).

\section{An application to functional equations}

In this section we give martingale characterization of the general measurable solution of the multidimensional Cauchy functional Equation.

{\bf Theorem 3.}
For a function $f: R^n\to R^m$ the following assertions are equivalent:

i) $f(x)=(f_1(x),...,f_m(x)), x\in R^n$ is a measurable function satisfying the  Cauchy functional equation
\begin{equation}\label{cosh}
f(x+y)=f(x)+f(y),
\end{equation}
for almost all $(x,y)$ in the sense of the Lebesgue measure on $R^{2n}$.

ii) $f(x)=(f_1(x),...,f_m(x)), x\in R^n$ is a measurable function such that the transformed process $(f(W_t), t\ge0)$ is a martingale satisfying condition C (of Theorem 1).

iii) $f(x)=Ax$    for some $m\times n$ constant matrix $A$ for almost all $x$ in the sense of the Lebesgue measure on $R^{n}$.

{\it Proof}. $i)\to ii)$. It is sufficient to show that $f_1(W)$ is a martingale satisfying condition $C)$.

Since for  each pair $(\xi,\eta)$ of random vectors with non-degenerate normal distribution  (with
density function $\rho(x,y)$)
$$
P\big(f_1(\xi+\eta)-f_1(\xi)-f_1(\eta)\neq 0\big)=
\int_{R^n}\int_{R^n} I_{(f_1(x+y)-f_1(x)-f_1(y)\neq 0)}\rho(x,y)dxdy=0,
$$
we have that
\begin{equation}\label{cosh}
f_1(\xi+\eta)=f_1(\xi)+f_1(\eta)\;\;\; a.s.
\end{equation}
Let $\tilde\xi,\tilde\eta$ be i.i.d. random vectors with normal distribution. Then, since the pairs $(\tilde\xi,\tilde\eta)$
and  $(\tilde\xi-\tilde\eta,\tilde\eta)$ have non-degenerate normal distributions, it follows from (\ref{cosh}) that
\begin{equation}\label{cosh1}
f_1(\tilde\xi+\tilde\eta)=f_1(\tilde\xi)+f_1(\tilde\eta)\;\;\; a.s.
\end{equation}
\begin{equation}\label{cosh2}
f_1(\tilde\xi-\tilde\eta)=f_1(\tilde\xi)-f_1(\tilde\eta)\;\;\; a.s.
\end{equation}
To show that $f_1(\tilde\xi)$ is integrable  we shall use the idea from \cite{SM} on application of the Bernstein theorem.

Let
$$
X=f_1(\tilde\xi)\;\;\;\text{and}\;\;\;Y=f_1(\tilde\eta).
$$
Then from (\ref{cosh1}) and (\ref{cosh2}) we have that
\begin{equation}\label{cosh3}
X+Y= f_1(\tilde\xi+\tilde\eta)\;\;\; a.s.
\end{equation}
\begin{equation}\label{cosh4}
X-Y=f_1(\tilde\xi-\tilde\eta)\;\;\; a.s.
\end{equation}
Since $\tilde\xi+\tilde\eta$ and $\tilde\xi-\tilde\eta$ are independent, the random variables $ f_1(\tilde\xi+\tilde\eta)$ and $f_1(\tilde\xi-\tilde\eta)$ will be also independent.
Therefore Bernstein's
theorem \cite{BE}(see also \cite{Q}) implies that $X=f_1(\tilde\xi)$ (and $Y=f_1(\tilde\eta))$ is normally distributed.
Hence the random variable $f_1(\tilde\xi)$ is square integrable.
In particular, this implies that
\begin{equation}\label{cosh5}
Ef_1^2(W_t)<\infty,\;\;\;\text{for every}\;\;\; t\ge0.
\end{equation}
Now let us show that  $Ef_1(\tilde\xi)=0$
for every  normally distributed random variable $\tilde\xi$ with zero mean.
Since $\tilde\xi$ and $\tilde\eta$ are independent with equal normal distribution we have that the random variables
$\tilde\xi+\tilde\eta, \tilde\xi-\tilde\eta$ and $\sqrt{2}\tilde\xi$ have the same normal distribution, hence
$$
Ef_1(\tilde\xi+\tilde\eta)=Ef_1 (\tilde\xi-\tilde\eta)=Ef_1 (\sqrt{2}\tilde\xi).
$$
Therefore,  from (\ref{cosh1}) and (\ref{cosh2}) we obtain that
$$
0=Ef_1 (\tilde\xi)-E f_1(\tilde\eta)=Ef_1 (\sqrt{2}\xi)=Ef_1 (\tilde\xi)+Ef_1 (\tilde\eta)=2Ef_1 (\tilde\xi).
$$
In particular,
\begin{equation}\label{cosh6}
Ef_1(W_t-W_s)=0\;\;\;\;\text{for all}\;\;\; s\le t.
\end{equation}

Substituting $\xi=W_t-W_s$ and $\eta=W_s$ in (\ref{cosh}) we have that
$$
f_1(W_t)-f_1(W_s)=f_1(W_t-W_s)\;\;\;a.s.
$$
Since $W_t-W_s$ is independent of ${\cal F}_s$, taking conditional expectations in this equality from (\ref{cosh6}) we obtain the
martingale equality
\beaa
E(f_1(W_t)-f_1(W_s)|{\cal F}_s)
=E(f_1(W_t-W_s))|{\cal F}_s)
=Ef_1(W_t-W_s)=0\;\;\;a.s
\eeaa
Denote $C(t)={\rm sign} (t)Ef_1^2(W_{|t|}),\;t\in R$. It follows from (\ref{cosh}) and the equality  $Ef_1(W_{t+s}-W_t)=0, \;t,s>0$ that
\beaa
C(t+s)=Ef_1^2(W_{t+s}-W_t+W_t)=Ef_1^2(W_{t+s}-W_t)\\
+2 Ef_1(W_{t+s}-W_t)Ef_1(W_t)+Ef_1^2(W_t)=C(t)+C(s),\;t,s>0.
\eeaa
From $C(t)=C(s)+C(t-s), \;t>s>0$ follows that
\beaa
C(t+(-s))=C(t-s)=C(t)-C(s)=C(t)+C(-s), \\
C(-t+s)=-C(t+(-s))=-C(t)-C(-s)=C(-t)+C(s),
\eeaa
which means $C(t+s)=C(t)+C(s)$ holds for all $t,s\in R.$
Thus, $C(t)$ is a bounded from below (on $R_+$) solution of the Cauchy one dimensional  functional equation
$$
C(t+s)=C(t)+C(s),\;\;\;\;s,t\in R.
$$
By well known result \cite{Died} the general solution bounded from below on some interval   is of the form  $C(t)=ct$ for some $c\in R$ .
Therefore, $Ef_1^2(W_t)=ct$ for some $c>0$.
From the equality $E(f_1(W_t)-f_1(W_s))^2|{\cal F}_s)=E(f_1(W_t)-f_1(W_s))^2=c(t-s)$  we obtain that $ct$ is the square characteristic of
the continuous modification of the martingale $f_1(W_t)$,
which means that conditions of Theorem 1 are satisfied.

$ii)\to iii)$. It follows from Theorem 1  follows that $f(x)=Ax + B$   for almost all $x$ in the sense of Lebesgue measure on $R^n$, for some $m\times n$ matrix $A$ and $B\in R^m$. This implies that
$$
f (\xi)=A\xi+B\;\;\;\;a.s.
$$
for every  normally distributed random variable $\xi$ with zero mean. Since
$Ef(\xi)=AE\xi=0$, taking mathematical expectations in the last equality we obtain that $B=0$. Thus,
$f(x)=Ax$   for almost all $x$ in the sense of Lebesgue measure.

$iii)\to i)$. If $f(x)=Ax\;a.e.$  for some $m\times n$ constant matrix $A$, for  independent random vectors $\xi$ and $\eta$
 having the standard normal distribution we have that
 $$
 P(f(\xi)=A\xi)=P(f(\eta)=A\eta)=P(f(\xi+\eta)=A(\xi+\eta))=1.
 $$
Therefore, avoiding three null sets we obtain from here
that $P-a.s.$
$$
f(\xi+\eta)=f(\xi)+ f(\eta),
$$
which implies  that (\ref{cosh}) is satisfied
for almost all $(x,y)$ in the sense of the Lebesgue measure on $R^{2n}$.

{\bf Remark.} In one-dimensional case Theorem 3 (more exactly, the equivalence of i) and iii)) follows from \cite{J} and \cite{DB}, where general additive functions are considered.

The following assertion is well known (see, e.g., \cite {AD}, \cite{SK}).  Theorem 3 gives a probabilistic proof of this result.

{\bf Corollary 3.}
 Let $f(x)=(f_1(x),...,f_m(x)), x\in R^n$ be a measurable function satisfying the  Cauchy functional equation
\begin{equation}\label{cosh3}
f(x+y)=f(x)+f(y),
\end{equation}
for  all $(x,y)$  $x, y\in R^{n}$.
Then $f(x)=Ax$  for some $m\times n$ constant matrix $A$.

{\it Proof}. It follows from the proof of Theorem 3 that the process $f(W_t)$ is a martingale which satisfies conditions of Theorem 1. From (\ref{cosh3}) we have that
$$
f(x+W_t)=f(x)+f(W_t),
$$
which implies that the process $f(x+W_t)$ is a martingale for any $x\in R^n$.
By  the martingale equality  we have that
$$
f(x)=E f(x+W_t)=\int_R f(x+y) \frac{1}{\sqrt{2\pi t}}e^{-\frac{y^2}{2t}}dy=
$$
 \begin{equation}\label{cont}
=\int_R f(y) \frac{1}{\sqrt{2\pi t}}e^{-\frac{(y-x)^2}{2t}}dy.
\end{equation}
Since $E|f(x+W_t)|<\infty$, equality (\ref{cont}) implies that the function  $f(x)$ is continuous and the proof follows from Corollary 1 of Theorem 1.

\appendix

\section{Appendix}

{\bf Lemma  A1}. Let $(X_t,t\ge 0)$ be a continuous function such that
\beaa
X_0=0,\;\limsup_{t\to\infty}X_t=\infty,\;
\liminf_{t\to\infty}X_t=-\infty.
\eeaa
If the composition $Y_t=h(X_t)$  is a right-continuous function, then the function $(h(x),x\in R)$ will be continuous.

 {\it Proof}.  It is sufficient to show, that $h$ is continuous on each  intervals $[0,b),(-b,0], b>0$. Let
 \beaa
 \tau(b)=\inf\{t>0; X_t\ge b\},\;\tau(b,0)=\inf\{t>\tau(b); X_t\le 0\},\\
 \a(x)=\sup\{t\le\tau(b); X_t= x\},\;\beta(x)=\sup\{\tau(b)<t<\tau(b,0); X_t= x\},\\
 0\le x<b.
 \eeaa
 It is evident, that $\a(x)$ is non-decreasing and right-continuous, and $\b(x)$  -- non-increasing and left-continuous. For example, the right-continuity of $a(x)$ is valid, since
if $x_n\downarrow x,\;\a(x_n)\downarrow\a^*\ge \a(x)$ then  $x_n=Y_{\a(x_n)}\downarrow Y_{\a^*}=x,\;\a^*\le \a(x).$
Equality $h(x)=Y_{\a(x)}=Y_{\b(x)}$ gives the continuity  of $h$  on $[0,b)$. Similarly can be shown the continuity
on $(-b,0]$.\qed

In particular, Lemma A1 implies that, if for a linear Brownian Motion $W$ the process $f(W_t)$ is a.s. right-continuous martingale, then the function $f=(f(x), x\in R)$ is continuous.
The following counterexample shows that in multidimensional case this fact is no longer true.

{\bf Counterexample.} Let $W=(W^1, W^2)$ be a two-dimensional Brownian Motion and let $h(x_1,x_2)=1_{(1,1)}(x_1,x_2)$.
Since
\begin{equation}\label{zero}
P\big(\omega: h(W_t)=0, \forall t\ge0\big)=P\big(\omega: W_t\neq (1,1), \forall t\ge 0\big)=0,
\end{equation}
by noting that the filtration ${\cal F}_t$ is complete, $(h(W_t), {\cal F}_t)$ is a martingale (the process indistinguishable from zero), but
the function $h$ is not continuous.
Note that, $h(x)=0$ almost everywhere with respect to the Lebesgue measure, which is in accordance with Theorem 1.

{\bf Remark.} If $h(x_1,x_2)= x_1I_{(R^2  - (1,1))}(x_1,x_2)$ then it follows from (\ref{zero}) that the process $f(W_t)$ is indistinguishable from the Brownian Motion $W^1_t$, but the function
$h$ is not continuous. Therefore, in theorem 2 from \cite{T} or the continuity of $f$ should be required, or the conclusion on almost surely linearity should be made.

In multidimensional case similar to Lemma A1 assertion will be valid, if we require the continuity of the composite function $h(X_t)$ for any
continuous function $X_t$ (and not only for almost all continuous paths with respect to the Wiener measure).

{\bf Lemma A2}. Let $h(x),x\in R^n$ be real-valued function such that $h(X_t)$ is continuous for each continuous $X_t,t\ge 0$. Then
$h$ is continuous.

{\it Proof}. Let $x_k,k=1,2,..$ be a convergent sequence and $x_0=\lim_{k\to\infty}x_k$.
The function defined by
$$X_t=\be{cases} x_0,\; {\rm if}\; t=0\\ x_{k+1}+k(k+1)(t-\frac1{k+1})(x_{k}-x_{k+1}),\;{\rm if}\;\frac1{k+1}<t\le\frac1k\\
x_1,\;{\rm if}\;t>1.
\ee{cases}$$
is continuous. Hence $h(X_t)$ is continuous and
$$
\lim_{k\to\infty}h(x_k)=\lim_{t\to0}h(X_t)=h(X_0)=h(x_0).$$

\end{document}